\documentclass[a4 paper,11pt,bezier]{article}
\usepackage[top=3cm,bottom=3cm,left=2.5cm,right=2.5cm]{geometry}
\usepackage{amsmath}
\usepackage{amsfonts,amsthm,amssymb}
\usepackage{amsfonts}
\usepackage{graphics,subfigure,caption2}
\usepackage{tikz, color}
\usepackage{graphics,epstopdf}
\usepackage{latexsym,bm}
\usepackage{amsfonts,amsthm,amssymb,mathrsfs,bbding}
\usepackage{CJK}

\usepackage{indentfirst}
\newtheorem{lem}{Lemma}[section]
\newtheorem{thm}[lem]{Theorem}
\newtheorem{cor}[lem]{Corollary}

\theoremstyle{plain}

\begin{document}

\title{The Restricted Edge-Connectivity of Strong Product Graphs\footnote{The research is supported by National Natural Science Foundation of China (12261086).}}
\author{Hazhe Ye, Yingzhi Tian\footnote{Corresponding author. E-mail: tianyzhxj@163.com (Y. Tian).}\\
{\small College of Mathematics and System Sciences, Xinjiang
University, Urumqi, Xinjiang, 830046, China}}
\date{}
\maketitle
		
\maketitle {\flushleft\bf Abstract:} The restricted edge-connectivity of a connected graph $G$, denoted by $\lambda^{\prime}(G)$, if it exists, is the minimum cardinality of a set of edges  whose deletion makes $G$ disconnected and each component with at least 2 vertices. It was proved that if $G$ is not a star and $|V(G)|\geq4$, then $\lambda^{\prime}(G)$ exists and $\lambda^{\prime}(G)\leq\xi(G)$, where $\xi(G)$ is the minimum edge-degree of $G$. Thus a graph $G$ is called maximally restricted edge-connected if $\lambda^{\prime}(G)=\xi(G)$; and a graph $G$  is called super restricted edge-connected if each minimum restricted edge-cut isolates an edge of $G$. 
The strong product of graphs $G$ and $H$, denoted by $G\boxtimes H$, is the graph with vertex set $V(G)\times V(H)$ and edge set $\{(x_1,y_1)(x_2,y_2)\ |\ x_1=x_2$ and $y_1y_2\in E(H)$; or $y_1=y_2$ and $x_1x_2\in E(G)$; or $x_1x_2\in E(G)$ and $y_1y_2\in E(H)$\}. 
In this paper, we determine, for any nontrivial connected graph $G$, the restricted edge-connectivity of $G\boxtimes P_n$, $G\boxtimes C_n$ and $G\boxtimes K_n$, where $P_n$, $C_n$ and $K_n$ are the path, the cycle and the complete graph on $n$ vertices, respectively. As corollaries, we give sufficient conditions for these  strong product graphs $G\boxtimes P_n$, $G\boxtimes C_n$ and $G\boxtimes K_n$ to be maximally restricted edge-connected and super restricted edge-connected.

\maketitle {\flushleft\textit{\bf Keywords}:} Edge-connectivity;  Restricted edge-connectivity; Strong product graphs; Maximally restricted edge-connected graphs;  Super restricted edge-connected graphs

\section{Introduction}

For notations and graph-theoretical terminology not defined here, we follow \cite{Bondy}. All graphs in this paper are undirected, simple and finite. Let $G=(V,E)$ be a graph, where $V=V(G)$ is the vertex set and $E=E(G)$ is the edge set.  The order of $G$ is $|V(G)|$, and the  size of $G$ is $e(G)=|E(G)|$.
For a vertex $u\in V(G)$, the neighborhood of $u$ in $G$ is $N_G(u)=\{v\in V(G)\ |\ v$ is adjacent to $u\}$, and the degree of $u$ in $G$ is $d_G(u)=|N_G(u)|$. The minimum degree $\delta(G)$ of $G$ is $min\{d_G(u)\ |\ u\in V(G)\}$. For an edge $e=uv\in E(G)$, $\xi_G(e)=d_G(u)+d_G(v)-2$ is the edge-degree of $e$ in $G$. The minimum edge-degree of $G$, denoted by $\xi(G)$, is $min\{d_G(u)+d_G(v)-2\ |\ e=uv\in E(G)\}$. Obviously, $\xi(G)\geq2\delta(G)-2$, with the equality holding if and only if there is an edge $e=uv\in E(G)$ such that $d_G(u)=d_G(v)=\delta(G)$. For a vertex set $A\subseteq V(G)$, the $induced$ $subgraph$ of $A$ in $G$, denoted by $G[A]$,  is the graph with vertex set $A$ and two vertices $u$ and $v$ in $A$ are adjacent if and only if they are adjacent in $G$.

For two nonempty subsets $X,Y\subseteq V(G)$, $[X,Y]_G$ denotes the set of edges with one end in $X$ and the other in $Y$. When $Y=V(G)\backslash X$, the set $[X,Y]_G$ is called an edge-cut of $G$ associated
with $X$. The edge-connectivity $\lambda(G)$ of a  graph $G$ is defined as the cardinality of a minimum edge-cut of $G$. It is well known that $\lambda(G)\leq\delta(G)$. Thus a graph $G$ is said to be maximally edge-connected if $\lambda(G)=\delta(G)$; and a graph $G$ is said to be super edge-connected if each minimum edge-cut isolates a vertex of $G$. A super edge-connected graph must be maximally edge-connected. But the converse is not true. For example, the cycle $C_n (n\geq4)$ is maximally edge-connected but not super edge-connected.

As an interconnection network can be modeled by a graph, the edge-connectivity can be used to measure the network reliability. But there is a deficiency, which allows all edges incident with a vertex to fail simultaneously. This situation is highly improbable in practical network applications.   
For compensating this deficiency, Esfahanian and Hakimi \cite{Esfahanian} introduced the notion of
restricted edge-connectivity. If an edge set $S \subseteq E(G)$ satisfies $G-S$ is disconnected and each component of $G-S$ has at least $2$ vertices, then $S$ is called a restricted edge-cut. If $G$ has at least one restricted edge-cut, then the restricted edge-connectivity of $G$, denoted by $\lambda^{\prime}(G)$, is the  cardinality of a minimum restricted edge-cut of $G$. It was proved in \cite{Esfahanian} that  if $G$ is not a star and its order is at least four, then $\lambda^{\prime}(G)\leq\xi(G)$.  Thus, if $\lambda^{\prime}(G)=\xi(G)$, then $G$ is said to be maximally restricted edge-connected; if each minimum restricted edge-cut isolates an edge of $G$, then $G$ is said to be super restricted edge-connected. A super restricted edge-connected graph must be maximally restricted edge-connected. But the converse is not true. For example, the cycle $C_n (n\geq6)$ is maximally restricted edge-connected but not super restricted edge-connected. 

The concept of graph product is utilized to construct larger graphs from smaller ones. There exist various kinds of graph products, including Cartesian product, direct product and  strong product, etc.  Given two graphs $G$ and $H$, the vertex sets of the Cartesian product $G\Box H$, the direct product $G\times H$ and the strong product $G\boxtimes H$ are all $V(G)\times V(H)$. For two  distinct  vertices $(x_1,y_1)$ and $(x_2,y_2)$, they  are adjacent in  $G\Box H$ if and only if $x_1=x_2$ and $y_1y_2\in E(H)$, or $y_1=y_2$ and $x_1x_2\in E(G)$; they are adjacent in $G\times H$ if and only if $x_1x_2\in E(G)$ and $y_1y_2\in E(H)$; and they are adjacent in $G\boxtimes H$ if and only if $x_1=x_2$ and $y_1y_2\in E(H)$, or $y_1=y_2$ and $x_1x_2\in E(G)$, or $x_1x_2\in E(G)$ and $y_1y_2\in E(H)$. Clearly, $E(G\boxtimes H)=E(G\Box H)\cup E(G\times H)$.

In \cite{Klavzar}, Klav\v{z}ar and \v{S}pacapan determined the edge-connectivity of the Cartesian product of two nontrivial graphs.  Shieh \cite{Shieh} characterized the super edge-connected Cartesian product graphs of two maximally edge-connected regular graphs. For the results on  the restricted edge-connectivity of Cartesian product graphs, see \cite{Liu,Lu,Ou1} for references.

Some bounds on the edge-connectivity of the direct product of graphs were given by Bre\v{s}ar and  \v{S}pacapan \cite{Bresar2}. The edge-connectivity of the direct product of a nontrivial graph and a complete graph was obtained by Cao,  Brglez,  \v{S}pacapan and  Vumar \cite{Cao}. In \cite{Spacapan}, \v{S}pacapan not only determined the edge-connectivity of the direct product of two general  graphs, but also characterized the structure of each minimum edge-cut in these direct product graphs. In \cite{Ma},  Ma, Wang and Zhang studied the restricted edge-connectivity of  the direct product of a nontrivial graph with a complete graph. In \cite{Bai}, Bai, Tian and Yin further studied the super restricted edge-connectedness of these direct product graphs.

In \cite{Bresar1}, Bre{\v{s}}ar and {\v{S}}pacapan determined the edge-connectivity of the strong products of two connected graphs. Ou and Zhao \cite{Ou2} studied the restricted edge-connectivity of strong product of two triangle-free graphs. In \cite{Wang}, Wang, Mao, Ye and Zhao gave an expression of the restricted edge-connectivity of the strong product graphs with two maximally restricted edge-connected graphs.

Motivated by the results above, we will study the restricted edge-connectivity of the strong product
of a nontrivial connected graph with a path, or a cycle, or a complete graph in this paper. As corollaries, we give sufficient conditions for these  strong product graphs to be maximally restricted edge-connected and super restricted edge-connected. In the next section, we will introduce some definitions and lemmas. The main results will be presented in Section 3.

\section{Preliminary}

Denote by $P_n$, $C_n$ and $K_n$  the path, the cycle and the complete graph on $n$ vertices, respectively.

Let $G$ and $H$ be two graphs. Define a natural projection $p$ on $V(G)\times V(H)$ as follows: $p(x,y)=y$ for any $(x,y)\in V(G)\times V(H)$. For any given $x\in V(G)$, the subgraph induced by $\{(x,y)|y\in V(H)\}$ in $G\boxtimes H$, denoted by $H^x$. Analogously, for any given $y\in V(H)$, the subgraph induced by $\{(x,y)|x\in V(G)\}$ in $G\boxtimes H$, denoted by $G^y$. Obviously, $H^x\cong H$ and $G^y\cong G$.

The edge-connectivity of the strong product of two connected nontrivial graphs was given in the following lemma.

\begin{lem}\label{Lemma 2.1}(\cite{Bresar1})
Let $G$ and $H$ be two connected nontrivial graphs. Then
$$
\lambda(G \boxtimes H)=\min \{\lambda(G)(|V(H)|+2e(H)), \lambda(H)(|V(G)|+2e(G)), \delta(G)+\delta(H)+\delta(G)\delta(H)\}.
$$
\end{lem}

Let $H$ be a connected graph. Define $K_2 \odot H=K_2 \boxtimes H-E(\{a\} \boxtimes H)-E(\{b\} \boxtimes H)$, where $V\left(K_2\right)=\{a, b\}$, $\{a\} \boxtimes H$ is the strong product of the complete graph with only one vertex $a$ and $H$, and $\{b\} \boxtimes H$ is the strong product of the complete graph with only one vertex $b$ and $H$. It is not difficult to see that $K_2 \odot H$ is connected if and only if $H$ is connected.

\begin{lem}\label{Lemma 2.2}(\cite{Yang})
Let $H$ be a connected graph and $S$ be an edge cut of $K_2 \odot H$, where $V\left(K_2\right)=\{a, b\}$. If the vertices of $\{a\}\boxtimes H$ are in different components of $K_2 \odot H-S$ as well as $\{b\} \boxtimes H$, then $|S| \geq 2 \lambda(H)$.
\end{lem}

Since for any $x\in X$, we have $\delta(G)\leq d_G(x)\leq |X|-1+|[X,\bar{X}]_G|$. Thus the following lemma holds.

\begin{lem}\label{Lemma 2.3}(\cite{Ou1})
Let $G$ be a connected graph. If $X$ is a nonempty subset of $V(G)$, then for any $x\in X$, we have  $|X|+|[X, \bar{X}]_G| \geq d_G(x)+1\geq \delta(G)+1$, with all equalities holding if and only if $X$  is a minimum-degree vertex.
\end{lem}

\begin{lem}\label{Lemma 2.4}(\cite{Esfahanian})
For any graph $G$ with order at least four and is not a star, then $\lambda^{\prime}(G)$ exists and $\lambda^{\prime}(G)\leq\xi(G)$.
\end{lem}

\begin{lem}\label{Lemma 2.5}(\cite{Ou2})
Let $G$ and $H$ be two connected graphs. Then
$$
\xi(G\boxtimes H)=\min \{\xi(G) \delta(H)+4 \delta(H)+\xi(G), \delta(G) \xi(H)+4 \delta(G)+\xi(H)\}.
$$
\end{lem}

\section{Main results}

\begin{thm}\label{theorem 3.1}
Let $G$ be a connected nontrivial graph with $m$ vertices. Then $\lambda^{\prime}(G\boxtimes P_n)$=min$\{(3n-2)\lambda(G), m+2e(G), 2\xi(G)+4, 5\delta(G)+1\}$, where $n\geq2$.
\end{thm}

\renewcommand\proofname{\bf{Proof}}

\begin{proof} Denote $\mathcal{G}=G\boxtimes P_n$. Let $V(G)=\{x_1, x_2, \ldots, x_m\}$ and $V(P_n)=\{y_1, y_2, \ldots, y_n\}$, where $y_jy_{j+1}\in E(P_n)$ for $j=1,2,\ldots,n-1$. Let $[X, \bar{X}]_G$ be a minimum edge-cut of $G$. Then $[X\times V(P_n), \bar{X}\times V(P_n)]_\mathcal{G}$ is a restricted edge-cut of $\mathcal{G}$. By $ |[X\times V(P_n), \bar{X}\times V(P_n)]_\mathcal{G}|=(3n-2)\lambda(G)$, we have $\lambda^{\prime}(\mathcal{G})\leq(3n-2)\lambda(G)$. Analogously, since $[\{y_1\},V(P_n)\backslash \{y_1\}]_{P_n}$ is a minimum edge-cut of $P_n$ and  $[\{y_1\}\times V(G), \{V(P_n)\backslash \{y_1\}\}\times V(G)]_\mathcal{G}$ is a restricted edge-cut of $\mathcal{G}$, we have $\lambda^{\prime}(\mathcal{G})\leq|[\{y_1\}\times V(G), \{V(P_n)\backslash \{y_1\}\}\times V(G)]_\mathcal{G}|=m+2e(G)$.
By Lemmas 2.4 and 2.5, $\lambda^{\prime}(\mathcal{G})\leq\xi(\mathcal{G})=\min \{\xi(G) \delta(P_n)+4 \delta(P_n)+\xi(G), \delta(G) \xi(P_n)+4 \delta(G)+\xi(P_n)\}=min\{2\xi(G)+4,5\delta(G)+1\}$.
Therefore, 
$\lambda^{\prime}(\mathcal{G})\leq min\{(3n-2)\lambda(G), m+2e(G),
2\xi(G)+4,5\delta(G)+1\}$.

Now, it is sufficient to prove $\lambda^{\prime}(\mathcal{G})\geq min\{(3n-2)\lambda(G), m+2e(G),
2\xi(G)+4,5\delta(G)+1\}$. Let $S$ be a minimum restricted edge-cut of $\mathcal{G}$. Then $\mathcal{G}-S$ has exactly two components, say $D_1$ and $D_2$, where $|V(D_1)|\geq2$ and $|V(D_2)|\geq2$. We consider two cases in the following.

\noindent{\bf Case 1.} Each vertex $x_i\in V(G)$  satisfies $P_n^{x_i}\cap D_1\neq\emptyset$ and $P_n^{x_i}\cap D_2\neq\emptyset$, or each vertex $y_j\in V(P_n)$  satisfies $G^{y_j}\cap D_1\neq\emptyset$ and $G^{y_j}\cap D_2\neq\emptyset$.

Assume each vertex $x_i\in V(G)$ satisfies $P_n^{x_i}\cap D_1\neq\emptyset$ and $P_n^{x_i}\cap D_2\neq\emptyset$. For $1\leq i\leq m$, denote $Y_i=V(P_n^{x_i})\cap V(D_1)$, $\bar{Y_i}=V(P_n^{x_i})\backslash Y_i$. By lemma 2.2, we have $|E(G[\{u,v\}]\odot P_n)\cap S|\geq2\lambda(P_n)$ for any edge $uv\in E(G)$. Thus

$$
\begin{aligned}
|S|&\geq \sum_{i=1}^m|[Y_i,\bar{Y_i}]_{P_n^{x_i}}|+\sum_{e=uv\in E(G)}|E(G[\{u,v\}]\odot P_n)\cap S|\\
&\geq m\cdot\lambda(P_n)+e(G)\cdot2\lambda(P_n)\\
&=m+2e(G).
\end{aligned}
$$

Analogously, if each vertex $y_j\in V(P_n)$ satisfies $G^{y_j}\cap D_1\neq\emptyset$ and $G^{y_j}\cap D_2\neq\emptyset$, then we have $|S|\geq(3n-2)\lambda(G)$.

\noindent{\bf Case 2.} There exists a vertex $x_a\in V(G)$ and a vertex $y_b\in V(P_n)$ such that $P_n^{x_a}\cap D_1=\emptyset$ and $G^{y_b}\cap D_1=\emptyset$, or $P_n^{x_a}\cap D_2=\emptyset$ and $G^{y_b}\cap D_2=\emptyset$.

Without loss of generality,  assume   $P_n^{x_a}\cap D_1=\emptyset$ and $G^{y_b}\cap D_1=\emptyset$. By the assumption, we know $V(P_n^{x_a})$ and $V(G^{y_b})$ are contained in $D_2$. Let $p(V(D_1))=\left\{y_{s+1}, y_{s+2},\ldots ,y_{s+k}\right\}$. Without loss of generality, assume $s+k<b$. For $1\leq i\leq k$, denote $X_i=V(G^{y_{s+i}})\cap V(D_1)$, $\bar{X_i}=V(G^{y_{s+i}})\backslash X_i$. For any $(x,y_{s+k})\in X_{k}$, we have $|[\{(x,y_{s+k})\},V(G^{y_{s+k+1}})]_\mathcal{G}|=d_G(x)+1$ by the definition of the strong product. Hence

\begin{align}
|S|&\geq\sum_{i=1}^k|[X_i,\bar{X_i}]_G{^{y_{s+i}}}|+|[X_k,V(G^{y_{s+k+1}})]_\mathcal{G}|+\sum_{i=1}^{k-1}|[X_i,
\bar{X}_{i+1}]_\mathcal{G}|+\sum_{i=1}^{k-1}|[X_{i+1},\bar{X_i}]_\mathcal{G}|\nonumber\\
&=\sum_{i=1}^k|[X_i,\bar{X_i}]_G{^{y_{s+i}}}|+\sum_{(x,y_{s+k})\in X_k}(d_G(x)+1)+\sum_{i=1}^{k-1}|[X_i,\bar{X}_{i+1}]_\mathcal{G}|+\sum_{i=1}^{k-1}|[X_{i+1},\bar{X_i}]_\mathcal{G}|.
\end{align}

\noindent{\bf Subcase 2.1.} $k=1$.

When $k=1$, we know that $V(D_1)\subseteq V(G^{y_{s+1}})$, that is, $V(D_1)=X_1$. Since $D_1$ is a connected graph, we have $G^{y_{s+1}}[X_1]$ is connected. If $|X_1|=2$, then $|S|\geq\xi(\mathcal{G})=min\{2\xi(G)+4,5\delta(G)+1\}$. If $|X_1|\geq3$, by $G^{y_{s+1}}[X_1]$ is connected, then there exists a vertex $(x,y_{s+1})\in X_1$ such that $d_G(x)\geq2$. Without loss of generality, assume   $(x_1,y_{s+1})\in X_1$ and $d_G(x_1)\geq2$. Let  $(x_2,y_{s+1}),(x_3,y_{s+1})\in N_{\mathcal{G}}((x_1,y_{s+1}))\cap X_1$. By lemma 2.3, $|[X_1,\bar{X_1}]_{G^{y_{s+1}}}|+|X_1|\geq d_{G}(x_1)+1$. Thus, by (1), we have

$$
\begin{aligned}
|S|&\geq|[X_1,\bar{X_1}]_{G^{y_{s+1}}}|+\sum_{(x,y_{s+1})\in X_1}(d_G(x)+1)\\
&=|[X_1,\bar{X_1}]_{G^{y_{s+1}}}|+|X_1|+\sum_{(x,y_{s+1})\in X_1}d_G(x)\\
&\geq d_{G}(x_1)+1+d_{G}(x_1)+d_{G}(x_2)+d_{G}(x_3)\\
&=(d_{G}(x_1)+d_{G}(x_2)-2)+( d_{G}(x_1)+d_{G}(x_3)-2)+5\\
&\geq2\xi(G)+5\\
&>2\xi(G)+4.
\end{aligned}
$$

\noindent{\bf Subcase 2.2.} $2\leq k\leq n-1$.

By Lemma 2.3 and $(1)$, we have
\begin{align}
|S|&\geq\sum_{i=1}^{k-1}|[X_i,\bar{X_i}]_{G^{y_{s+i}}}|+|[X_k,\bar{X_k}]_{G^{y_{s+k}}}|+
|X_k|(\delta(G)+1)+\sum_{i=1}^{k-1}|[X_i,\bar{X}_{i+1}]_\mathcal{G}|+
\sum_{i=1}^{k-1}|[X_{i+1},\bar{X_i}]_\mathcal{G}|\nonumber\\
&=\sum_{i=1}^{k-1}|[X_i,\bar{X_i}]_{G^{y_{s+i}}}|+|[X_k,\bar{X_k}]_{G^{y_{s+k}}}|+|X_k|
+|X_k|\delta(G)+\sum_{i=1}^{k-1}|[X_i,\bar{X}_{i+1}]_\mathcal{G}|+
\sum_{i=1}^{k-1}|[X_{i+1},\bar{X_i}]_\mathcal{G}|\nonumber\\
&\geq(k-1)\lambda(G)+\delta(G)+1+|X_k|\delta(G)+\sum_{i=1}^{k-1}|[X_i,\bar{X}_{i+1}]_\mathcal{G}|+
\sum_{i=1}^{k-1}|[X_{i+1},\bar{X_i}]_\mathcal{G}|.
\end{align}

For any $(x,y_{s+k-1})\in X_{k-1}$, we have $|[\{(x,y_{s+k-1})\},\bar{X_k}]_\mathcal{G}|\geq(\delta(G)+1)-|X_k|$, and for any $(x,y_{s+k})\in X_{k}$, we have $|[\{(x,y_{s+k})\},\bar{X}_{k-1}]_\mathcal{G}|\geq(\delta(G)+1)-|X_{k-1}|$ .
It follows that
\begin{align}
&\sum_{i=1}^{k-1}|[X_i,\bar{X}_{i+1}]_\mathcal{G}|+\sum_{i=1}^{k-1}|[X_{i+1},\bar{X_i}]_\mathcal{G}|\nonumber\\
&\geq|[X_{k-1},\bar{X_{k}}]_\mathcal{G}|+|[X_{k},\bar{X}_{k-1}]_\mathcal{G}|\nonumber\\
&\geq|X_{k-1}|[(\delta(G)+1)-|X_{k}|]+|X_{k}|[(\delta(G)+1)-|X_{k-1}|].
\end{align}

If $|X_k|\geq4$, then, by (2), we have $|S|>5\delta(G)+1$. So we just need to consider $1\leq|X_k|\leq3$. There are three subcases in the following.

\noindent{\bf Subcase 2.2.1.} $|X_k|=1$.

If $|X_{k-1}|=1$, then  $|[X_{k-1},\bar{X}_{k-1}]_\mathcal{G}|\geq\delta(G)$ and
$|[X_{k},\bar{X_{k}}]_\mathcal{G}|\geq\delta(G)$. 
Furthermore, by $(3)$,  we have $\sum_{i=1}^{k-1}|[X_i,\bar{X}_{i+1}]_\mathcal{G}|+\sum_{i=1}^{k-1}|[X_{i+1},\bar{X_i}]_\mathcal{G}|
\geq2\delta(G)$. Thus by $(1)$, we obtain that $|S|\geq\sum_{i=1}^{k-2}|[X_i,\bar{X_i}]_G{^{y_{s+i}}}|+|[X_{k-1},\bar{X}_{k-1}]_\mathcal{G}|
+|[X_{k},\bar{X_{k}}]_\mathcal{G}|+\sum_{(x,y_{s+k})\in X_k}(d_G(x)+1)+\sum_{i=1}^{k-1}|[X_i,\bar{X}_{i+1}]_\mathcal{G}|
+\sum_{i=1}^{k-1}|[X_{i+1},\bar{X_i}]_\mathcal{G}|
\geq(k-2)\lambda(G)+\delta(G)+\delta(G)+\delta(G)+1+2\delta(G)\geq5\delta(G)+1$.

If $|X_{k-1}|\geq2$, then by $(2)$ and $(3)$, we have
$|S|\geq(k-1)\lambda(G)+\delta(G)+1+\delta(G)+|X_{k-1}|[(\delta(G)+1)-1]+(\delta(G)+1)-|X_{k-1}|
=(k-1)\lambda(G)+|X_{k-1}|(\delta(G)-1)+3\delta(G)+2
\geq5\delta(G)+1.$

\noindent{\bf Subcase 2.2.2.} $|X_k|=2$.

If $\delta(G)=1$, then $\lambda(G)=\delta(G)=1$. By lemma 2.2, we obtain that $\sum_{i=1}^{k-1}|[X_i,\bar{X}_{i+1}]_\mathcal{G}|+\sum_{i=1}^{k-1}|[X_{i+1},
\bar{X_i}]_\mathcal{G}|\geq(k-1)\cdot2\lambda(G)=2(k-1)$. By $(1)$, we have
$|S|\geq k\lambda(G)+2(\delta(G)+1)+2(k-1)=3k+2\geq8>5\delta(G)+1$.

If $\delta(G)=2$, then by $(3)$, we have $|[X_{k-1},\bar{X_{k}}]_\mathcal{G}|+|[X_{k},\bar{X}_{k-1}]_\mathcal{G}|\geq
|X_{k-1}|[(\delta(G)+1)-|X_{k}|]+|X_{k}|[(\delta(G)+1)-|X_{k-1}|]
=-|X_{k-1}|+6\geq4$ when $|X_{k-1}|\leq2$. Suppose $|X_{k-1}|\geq3$. Since $|X_{k}|=2$ and $\delta(G)=2$, we have  $|[X_{k-1},\bar{X_{k}}]_\mathcal{G}|\geq |X_{k-1}|(\delta(G)+1-|X_{k}|)\geq3$. Thus, by lemma 2.2, we have
$\sum_{i=1}^{k-1}|[X_i,\bar{X}_{i+1}]_\mathcal{G}|+\sum_{i=1}^{k-1}|[X_{i+1},\bar{X_i}]_\mathcal{G}|
\geq\sum_{i=1}^{k-2}|[X_i,\bar{X}_{i+1}]_\mathcal{G}|+\sum_{i=1}^{k-2}|[X_{i+1},\bar{X_i}]_\mathcal{G}|+
|[X_{k-1},\bar{X_{k}}]_\mathcal{G}|+|[X_{k},\bar{X}_{k-1}]_\mathcal{G}|
\geq(k-2)\cdot2\lambda(G)+3
\geq2k-1.$
Hence, by $(1)$, we have
$|S|\geq k\lambda(G)+2(\delta(G)+1)+2k-1
=3k+5\geq11\geq5\delta(G)+1$.

If $\delta(G)\geq3$, then by $(2)$ and $(3)$, we have
$|S|\geq(k-1)\lambda(G)+\delta(G)+1+2\delta(G)+|X_{k-1}|[(\delta(G)+1)-2]+2[(\delta(G)+1)-|X_{k-1}|]
=(k-1)\lambda(G)+|X_{k-1}|(\delta(G)-3)+(5\delta(G)+3)
>5\delta(G)+1$.

\noindent{\bf Subcase 2.2.3}  $|X_k|=3$.

If $\delta(G)=1$, then $\lambda(G)=\delta(G)=1$. By lemma 2.2, we obtain that  $\sum_{i=1}^{k-1}|[X_i,\bar{X}_{i+1}]_\mathcal{G}|+\sum_{i=1}^{k-1}|[X_{i+1},\bar{X_i}]_\mathcal{G}|
\geq(k-1)\cdot2\lambda(G)=2(k-1)$. Thus, by $(1)$, we have $|S|\geq k\lambda(G)+3(\delta(G)+1)+2(k-1)=3k+4\geq10>5\delta(G)+1$.

Now we consider $\delta(G)\geq2$. Let $X_k=\{(x_1,y_{s+k}),(x_2,y_{s+k}),(x_3,y_{s+k})\}$. If $G^{y_{s+k}}[x_1,x_2,x_3]$ is connected, then, by a similar argument as Subcase 2.1, we can also obtain $|S|>2\xi(G)+4$. If $G^{y_{s+k}}[x_1,x_2,x_3]$ is not connected, then $G^{y_{s+k}}[x_1,x_2,x_3]$ must contain isolated vertices. Therefore, $|[X_k,\bar{X_k}]_{G^{y_{s+k}}}|\geq\delta(G)$. By $(3)$, it follows that  $\sum_{i=1}^{k-1}|[X_i,\bar{X}_{i+1}]_\mathcal{G}|+\sum_{i=1}^{k-1}|[X_{i+1},\bar{X_i}]_\mathcal{G}|\geq |[X_{k-1},\bar{X_{k}}]_\mathcal{G}|\geq|X_{k-1}|[(\delta(G)+1)-3]$.
Hence, by $(1)$, we obtain
$|S|\geq\sum_{i=1}^{k-1}|[X_i,\bar{X_i}]_{G^{y_{s+i}}}|+|[X_k,\bar{X_k}]_{G^{y_{s+k}}}|+
|X_k|(\delta(G)+1)+\sum_{i=1}^{k-1}|[X_i,\bar{X}_{i+1}]_\mathcal{G}|+
\sum_{i=1}^{k-1}|[X_{i+1},\bar{X_i}]_\mathcal{G}|
\geq(k-1)\lambda(G)+\delta(G)+3(\delta(G)+1)+|X_{k-1}|[(\delta(G)+1)-3]
=(k-1)\lambda(G)+(|X_{k-1}|-1)(\delta(G)-2)+5\delta(G)+1
>5\delta(G)+1.$

This proof is thus complete.
\end{proof}

Since $\xi(G\boxtimes P_n)=min\{2\xi(G)+4,5\delta(G)+1\}$, Theorem 3.1 implies the following corollary.

\begin{cor}\label{Corollary 3.2}
Let $G$ be connected nontrivial graph with $m$ vertices. If min$\{(3n-2)\lambda(G), m+2e(G)\}\geq min\{2\xi(G)+4,5\delta(G)+1\}$, then $G\boxtimes P_n$ is maximally restricted edge-connected, where $n\geq2$.
\end{cor}

By checking through the proof of Theorem 3.1, we find that $|S|\geq$min$\{(3n-2)\lambda(G), m+2e(G)\}$ or 
$|S|>$min$\{2\xi(G)+4,5\delta(G)+1\}$ when both $|V(D_1)|\geq3$ and $|V(D_2)|\geq3$. Thus we have the following corollary.

\begin{cor}\label{Corollary 3.3}
Let $G$ be a connected nontrivial graph with $m$ vertices. If min$\{(3n-2)\lambda(G), m+2e(G)\}>min\{2\xi(G)+4,5\delta(G)+1\}$, then $G\boxtimes P_n$ is super restricted edge-connected, where $n\geq2$.
\end{cor}

\begin{thm}\label{theorem 3.4}
Let $G$ be a connected nontrivial graph with $m$ vertices. Then $\lambda^{\prime}(G\boxtimes C_n)=min\{3n\lambda(G), 2(m+2e(G)), 6\delta(G)+2\}$, where $n\geq3$.
\end{thm}

\renewcommand\proofname{\bf{Proof}}

\begin{proof}Denote $\mathcal{G}=G\boxtimes C_n$. Let $V(G)=\{x_1, x_2, \ldots, x_m\}$ and $V(C_n)=\{y_1,y_2,\ldots,y_n\}$, where $y_jy_{j+1}\in E(C_n)$ for $j=1,2,\ldots,n$ ($y_{n+1}=y_1$). Let $[X, \bar{X}]_G$ be a minimum edge-cut of $G$. Then $[X\times V(C_n), \bar{X}\times V(C_n)]_\mathcal{G}$ is a restricted edge-cut of $\mathcal{G}$. By $|[X\times V(C_n), \bar{X}\times V(C_n)]_\mathcal{G}|=3n\lambda(G)$, we have $\lambda^{\prime}(\mathcal{G})\leq3n\lambda(G)$. Analogously, since $[\{y_1\}, V(C_n)\backslash \{y_1\}]_{C_n}$ is a minimum edge-cut of $C_n$, we have $\lambda^{\prime}(\mathcal{G})\leq|[\{y_1\}\times V(G), \{V(C_n)\backslash \{y_1\}\}\times V(G)]_\mathcal{G}|=2(m+2e(G))$ by $[\{y_1\}\times V(G), \{V(C_n)\backslash \{y_1\}\}\times V(G)]_\mathcal{G}$ is a restricted edge-cut of $\mathcal{G}$.
By Lemmas 2.4 and 2.5, 
$\lambda^{\prime}(\mathcal{G})\leq\xi(\mathcal{G})=min\{\xi(G)\delta(C_n)+4\delta(C_n)+\xi(G),  \delta(G)\xi(C_n)+4\delta(G)+\xi(C_n)\}
=min\{3\xi(G)+8, 6\delta(G)+2\}
=6\delta(G)+2$.
Therefore, $\lambda^{\prime}(\mathcal{G})\leq min\{3n\lambda(G), 2(m+2e(G)), 6\delta(G)+2\}$.

Now, it is sufficient to prove $\lambda^{\prime}(\mathcal{G})\geq min\{3n\lambda(G), 2(m+2e(G)), 6\delta(G)+2\}$. Let $S$ be a minimum restricted edge-cut of $\mathcal{G}$. Then $\mathcal{G}-S$ has exactly two components, say $D_1$ and $D_2$, where $|V(D_1)|\geq2$ and $|V(D_2)|\geq2$. We consider two cases in the following.

\noindent{\bf Case 1.} Each vertex $x_i\in V(G)$ satisfies $C_n^{x_i}\cap D_1\neq\emptyset$ and $C_n^{x_i}\cap D_2\neq\emptyset$, or each vertex $ y_j\in V(C_n)$ satisfies $G^{y_j}\cap D_1\neq\emptyset$ and $G^{y_j}\cap D_2\neq\emptyset$.

Assume each vertex $x_i\in V(G)$ satisfies $C_n^{x_i}\cap D_1\neq\emptyset$ and $C_n^{x_i}\cap D_2\neq\emptyset$. For $1\leq i\leq m$, denote $Y_i=V(C_n^{x_i})\cap V(D_1)$, $\bar{Y_i}=V(C_n^{x_i})\backslash Y_i$. By lemma 2.2, we have $|E(G[\{u,v\}]\odot C_n)\cap S|\geq2\lambda(C_n)$ for any edge $uv\in E(G)$. Thus

$$
\begin{aligned}
|S|&\geq \sum_{i=1}^m|[Y_i,\bar{Y_i}]_{C_n^{x_i}}|+\sum_{e=uv\in E(G)}|E(G[\{u,v\}]\odot C_n)\cap S|\\
&\geq m\cdot\lambda(C_n)+e(G)\cdot2\lambda(C_n)\\
&=2(m+2e(G)).
\end{aligned}
$$

Analogously, if each vertex $y_j\in V(C_n)$ satisfies $G^{y_j}\cap D_1\neq\emptyset$ and $G^{y_j}\cap D_2\neq\emptyset$, then we have $|S| \geq3n\lambda(G)$.

\noindent{\bf Case 2.} There exists a vertex $x_a\in V(G)$ and a vertex $y_b\in V(C_n)$ such that $C_n^{x_a}\cap D_1=\emptyset$ and $G^{y_b}\cap D_1=\emptyset$, or $C_n^{x_a}\cap D_2=\emptyset$ and $G^{y_b}\cap D_2=\emptyset$.

Without loss of generality,  assume  $C_n^{x_a}\cap D_1=\emptyset$ and $G^{y_b}\cap D_1=\emptyset$. By the assumption, we know $V(C_n^{x_a})$ and $V(G^{y_b})$ are contained in $D_2$. Let $p(V(D_1))=\left\{y_{s+1},y_{s+2},\ldots ,y_{s+k}\right\}$, where the addition is modular $n$ operation. Without loss of generality, assume $s+k<b$. For $1\leq i\leq k$, denote $X_i=V(G^{y_{s+i}})\cap V(D_1)$, $\bar{X_i}=V(G^{y_{s+i}})\backslash X_i$. For any vertex $(x,y_{s+1})\in X_1$, we have $|[\{(x,y_{s+1})\},V(G^{y_{s}})]_\mathcal{G}|= d_G(x)+1\geq\delta(G)+1$. Then $|[X_1,V(G^{y_{s}})]_\mathcal{G}|\geq|X_1|(\delta(G)+1)$. Analogously, $|[X_k,V(G^{y_{s+k+1}})]_\mathcal{G}|\geq|X_k|(\delta(G)+1)$. Hence

\begin{align}
|S|&\geq\sum_{i=1}^k|[X_i,\bar{X_i}]_G{^{y_{s+i}}}|+|[X_1,G^{y_s}]_\mathcal{G}|+|[X_k,G^{y_{s+k+1}}]_\mathcal{G}|+
\sum_{i=1}^{k-1}|[X_i,\bar{X}_{i+1}]_\mathcal{G}|+\sum_{i=1}^{k-1}|[X_{i+1},\bar{X_i}]_\mathcal{G}|\nonumber\\
&\geq\sum_{i=1}^{k}|[X_i,\bar{X_i}]_G{^{y_{s+i}}}|+|X_1|(\delta(G)+1)+|X_k|(\delta(G)+1)+
\sum_{i=1}^{k-1}|[X_i,\bar{X}_{i+1}]_\mathcal{G}|
+\sum_{i=1}^{k-1}|[X_{i+1},\bar{X_i}]_\mathcal{G}|.
\end{align}

\noindent{\bf Subcase 2.1.} $k=1$.

When $k=1$, we know that $V(D_1)\subseteq V(G^{y_{s+1}})$, that is, $V(D_1)=X_1$. Since $D_1$ is a connected graph, we have $G^{y_{s+1}}[X_1]$ is connected. If $|X_1|=2$, then $|S|\geq\xi(\mathcal{G})=6\delta(G)+2$. If $|X_1|\geq3$, then by lemma 2.3 and $(4)$, we have
$$
\begin{aligned}
|S|&\geq|[X_1,\bar{X_1}]_G{^{y_{s+1}}}|+|X_1|(\delta(G)+1)+|X_1|(\delta(G)+1)\\
&=|[X_1,\bar{X_1}]_G{^{y_{s+1}}}|+|X_1|+|X_1|\delta(G)+|X_1|(\delta(G)+1)\\
&\geq\delta(G)+1+3\delta(G)+3(\delta(G)+1)\\
&=7\delta(G)+4\\
&>6\delta(G)+2.
\end{aligned}
$$

\noindent{\bf Subcase 2.2.} $2\leq k\leq n-1$.

By Lemma 2.3 and $(4)$, we have

\begin{align}
|S|&\geq\sum_{i=2}^{k-1}|[X_i,\bar{X_i}]_G{^{y_{s+i}}}|+|[X_1,\bar{X_1}]_G{^{y_{s+1}}}|+|X_1|+
|[X_k,\bar{X_k}]_G{^{y_{s+k}}}|+|X_k|\nonumber\\
&+|X_1|\delta(G)+|X_k|\delta(G)+\sum_{i=1}^{k-1}|[X_i,\bar{X}_{i+1}]_\mathcal{G}|+
\sum_{i=1}^{k-1}|[X_{i+1},\bar{X_i}]_\mathcal{G}|\nonumber\\
&\geq(k-2)\lambda(G)+2\delta(G)+2+(|X_1|+|X_k|)\delta(G)+\sum_{i=1}^{k-1}|[X_i,\bar{X}_{i+1}]_\mathcal{G}|
+\sum_{i=1}^{k-1}|[X_{i+1},\bar{X_i}]_\mathcal{G}|
\end{align}

If $|X_1|+|X_k|\geq4$, then we have $|S|>6\delta(G)+2$. So we just need to consider $2\leq(|X_1|+|X_k|\leq3$, that is, $|X_1|=1$ and $|X_k|=1$, or $|X_1|=1$ and $|X_k|=2$, or $|X_1|=2$ and $|X_k|=1$. Without loss of generality, assume $|X_{k}|=1$. By $(3)$ and $(5)$, we have
$$
\begin{aligned}
|S|&\geq(k-2)\lambda+2\delta(G)+2+2\delta(G)+|X_{k-1}|[(\delta(G)+1)-1]+[(\delta(G)+1)-|X_{k-1}|]\\
&=(k-2)\lambda+(|X_{k-1}|-1)(\delta(G)-1)+(6\delta(G)+2)\\
&\geq6\delta(G)+2.
\end{aligned}
$$
This proof is thus complete.
\end{proof}

Since $\xi(G\boxtimes C_n)=6\delta(G)+2$, Theorem 3.4 implies the following corollary.

\begin{cor}\label{Corollary 3.5}
Let $G$ be a connected nontrivial graph with $m$ vertices. If min$\{3n\lambda, 2(m+2e(G))\}\geq 6\delta(G)+2$, then $G\boxtimes C_n$ is maximally restricted edge-connected, where $n\geq3$.
\end{cor}

By checking through the proof of Theorem 3.4, we find that $|S|\geq$min$\{3n\lambda, 2(m+2e(G))\}$ or 
$|S|>6\delta(G)+2$ when both $|V(D_1)|\geq3$ and $|V(D_2)|\geq3$. Thus we have the following corollary.

\begin{cor}\label{Corollary 3.6}
Let $G$ be a connected nontrivial graph with $m$ vertices. If min$\{3n\lambda, 2(m+2e(G))\}> 6\delta(G)+2$, then $G\boxtimes C_n$ is super restricted edge-connected, where $n\geq3$.
\end{cor}

\begin{thm}\label{theorem 3.7}
Let $G$ be a connected nontrivial graph with $m$ vertices. Then $\lambda^{\prime}(G\boxtimes K_n)=min\{n^2\lambda(G), (n-1)(m+2e(G)), 2n\delta(G)+2n-4\}$, where $n\geq4$.
\end{thm}

\renewcommand\proofname{\bf{Proof}}

\begin{proof} Denote $\mathcal{G}=G\boxtimes K_n$. Let $V(G)=\{x_1, x_2, \ldots, x_m\}$ and $V(K_n)=\{y_1, y_2, \ldots , y_n\}$. Let $[X, \bar{X}]_G$ be a minimum edge-cut of $G$. Then $[X\times V(K_n), \bar{X}\times V(K_n)]_\mathcal{G}$ is a restricted edge-cut of $\mathcal{G}$. By $ |[X\times V(K_n), \bar{X}\times V(K_n)]_\mathcal{G}|=n^2\lambda(G)$, we have $\lambda^{\prime}(\mathcal{G})\leq n^2\lambda(G)$. Analogously, since $[\{y_1\}, V(K_n)\backslash \{y_1\}]_{K_n}$ is a minimum edge-cut of $K_n$, we have
 $\lambda^{\prime}(\mathcal{G})\leq|[\!\{y_1\}\times V(G)\!,\!\{V(K_n)\backslash \{y_1\}\}\! \times V(G)]_\mathcal{G}|=(n-1)(m+2e(G))$
by $[\{y_1\}\times V(G), \!\{V(K_n)\backslash \{y_1\}\}\times V(G)\!]_\mathcal{G}$ is a restricted edge-cut of $\mathcal{G}$. 
By Lemmas 2.4 and 2.5, we have
$\lambda^{\prime}(\mathcal{G})\leq
\xi(\mathcal{G})=min\{\xi(G)\delta(K_n)+4\delta(K_n)+\xi(G), \delta(G)\xi(K_n)+4\delta(G)+\xi(K_n)\}
=min\{n\xi(G)+4n-4, 2n\delta(G)+2n-4\}
=2n\delta(G)+2n-4$.
Therefore,
$\lambda^{\prime}(\mathcal{G})\leq min\{n^2\lambda(G), (n-1)(m+2e(G)), 2n\delta(G)+2n-4\}$.

Now, it is sufficient to prove $\lambda^{\prime}(\mathcal{G})\geq min\{n^2\lambda, (n-1)(m+2e(G)), 2n\delta(G)+2n-4\}$. Let $S$ be a minimum restricted edge-cut of $\mathcal{G}$. Then $\mathcal{G}-S$ has exactly two components, say $D_1$ and $D_2$, where $|V(D_1)|\geq2$ and $|V(D_2)|\geq2$. We consider two cases in the following.

\noindent{\bf Case 1.} Each vertex $x_i\in V(G)$  satisfies $K_n^{x_i}\cap D_1\neq\emptyset$ and $K_n^{x_i}\cap D_2\neq\emptyset$, or each vertex $y_j\in V(K_n)$ satisfies $G^{y_j}\cap D_1\neq\emptyset$ and $G^{y_j}\cap D_2\neq\emptyset$.

Assume each vertex $x_i\in V(G)$ satisfies $K_n^{x_i}\cap D_1\neq\emptyset$ and $K_n^{x_i}\cap D_2\neq\emptyset$. For $1\leq i\leq m$, denote $Y_i=V(K_n^{x_i})\cap V(D_1)$, $\bar{Y_i}=V(K_n^{x_i})\backslash Y_i$. By lemma 2.2, we have $|E(G[\{u,v\}]\odot K_n)\cap S|\geq2\lambda(K_n)$ for any edge $uv\in E(G)$. Thus
$$
\begin{aligned}
|S|&\geq \sum_{i=1}^m|[Y_i,\bar{Y_i}]_{K_n^{x_i}}|+\sum_{e=uv\in E(G)}|E(G[\{u,v\}]\odot K_n)\cap S|\\
&\geq m\cdot\lambda(K_n)+e(G)\cdot2\lambda(K_n)\\
&=(n-1)(m+2e(G)).
\end{aligned}
$$

Analogously, if each vertex $y_j\in V(K_n)$ satisfies $G^{y_j}\cap D_1\neq\emptyset$ and $G^{y_j}\cap D_2\neq\emptyset$, then we have $|S| \geq n^2\lambda(G)$.

\noindent{\bf Case 2.} There exists a vertex $x_a\in V(G)$ and a vertex $y_b\in V(K_n)$ such that $K_n^{x_a}\cap D_1=\emptyset$ and $G^{y_b}\cap D_1=\emptyset$, or $K_n^{x_a}\cap D_2=\emptyset$ and $G^{y_b}\cap D_2=\emptyset$.

Without loss of generality, assume $K_n^{x_a}\cap D_1=\emptyset$ and $G^{y_b}\cap D_1=\emptyset$. By the assumption, we know $V(K_n^{x_a})$ and $V(G^{y_b})$ are contained in $D_2$. Since any two distinct vertices are adjacent in $K_n$, by renaming the vertices of $V(K_n)$,  we can let $p(V(D_1))=\left\{y_{s+1},y_{s+2},\ldots, y_{s+k}\right\}$. Furthermore, assume  $s+k<b$. For $1\leq i\leq k$, let $X_i=V(G^{y_{s+i}})\cap V(D_1)$, $\bar{X_i}=V(G^{y_{s+i}})\backslash X_i$. Denote $Y=V(K_n)\setminus p(V(D_1))$. By the definition of the strong product, for any $y\in Y$, we have $|[X_i,G^y]_\mathcal{G}|\geq|X_i|(\delta(G)+1)$. Hence

\begin{align}
|S|&\geq\sum_{i=1}^k|[X_i,\bar{X_i}]_{G^{y_{s+i}}}|+\sum_{i=1}^k\sum_{y\in Y}|[X_i,G^y]_\mathcal{G}|+
\sum_{i=1}^{k}\sum_{j\in\{1,\ldots,k\}\setminus\{i\}}|[X_i,\bar{X_{j}}]_\mathcal{G}|\nonumber\\
&\geq\sum_{i=1}^k|[X_i,\bar{X_i}]_{G^{y_{s+i}}}|+\sum_{i=1}^k|X_i|(\delta(G)+1)(n-k)
+\sum_{i=1}^{k}\sum_{j\in\{1,\ldots,k\}\setminus\{i\}}|[X_i,\bar{X_{j}}]_\mathcal{G}|\\
&=\sum_{i=1}^k|[X_i,\bar{X_i}]_{G^{y_{s+i}}}|+\sum_{i=1}^k|X_i|+\sum_{i=1}^k|X_i|[(\delta(G)+1)(n-k)-1]
+\sum_{i=1}^{k}\sum_{j\in\{1,\ldots,k\}\setminus\{i\}}|[X_i,\bar{X_{j}}]_\mathcal{G}|\nonumber\\
&\geq k(\delta(G)+1)+\sum_{i=1}^k|X_i|[(\delta(G)+1)(n-k)-1]
+\sum_{i=1}^{k}\sum_{j\in\{1,\ldots,k\}\setminus\{i\}}|[X_i,\bar{X_{j}}]_\mathcal{G}|.
\end{align}

\noindent{\bf Subcase 2.1.} $k=1$.

When $k=1$, we know that $V(D_1)\subseteq V(G^{y_{s+1}})$, that is, $V(D_1)=X_1$.  Since $D_1$ is a connected graph, we have $G^{y_{s+1}}[X_1]$ is connected. If $|X_1|=2$, then $|S|\geq\xi(\mathcal{G})=2n\delta(G)+2n-4$. If $|X_1|\geq3$, then by $(7)$, we have
$$
\begin{aligned}
|S|&\geq\delta(G)+1+3[(\delta(G)+1)(n-1)-1]\\
&=(n-2)\delta(G)+(n-1)+2n\delta(G)+2n-4\\
&>2n\delta(G)+2n-4.
\end{aligned}
$$

\noindent{\bf Subcase 2.2.} $2\leq k\leq n-1$.

\noindent{\bf Subcase 2.2.1.}   For each  $i\in\{1,\ldots,k\}$, $|X_i|\geq2$.

By lemma 2.2, we have
$$
\begin{aligned}
\sum_{i=1}^{k}\sum_{j\in\{1,\ldots,k\}\setminus\{i\}}|[X_i,\bar{X_{j}}]_\mathcal{G}|&\geq \frac{k(k-1)}{2}2\lambda(G)\geq k(k-1).
\end{aligned}
$$
Since $|X_i|\geq2$, we know that $\sum_{i=1}^k|X_i|\geq2k$.
By $(7)$, we have
$$
\begin{aligned}
|S|&\geq k(\delta(G)+1)+2k[(\delta(G)+1)(n-k)-1]+k(k-1)\\
&=k\delta(G)+k+2kn\delta(G)-2k^2\delta(G)+2kn-2k^2-2k+k^2-k\\
&=(-2k^2+2kn+k)\delta(G)+(-k^2+2kn-2k).
\end{aligned}
$$

Let $f_1(k)=-2k^2+2kn+k$ and $f_2(k)=-k^2+2kn-2k$. Since $2\leq k\leq n-1$ and $n\geq4$,
we have $f_1(k)\geq min\{f_1(2),f_1(n-1)\}=min\{4n-6,3n-3\}>2n$ and $f_2(k)\geq min\{f_2(2),f_2(n-1)\}=min\{4n-8,(n-2)^2+2n-3\}>2n-4$.  Thus we obtain  $|S|>2n\delta(G)+2n-4$.

\noindent{\bf Subcase 2.2.2.} There are at least two integers in $\{1,\ldots,k\}$, say $1$ and $2$, such that 
$|X_1|=|X_2|=1$.

By $|X_1|=1$, we have $|[X_i,\bar{X_1}]_\mathcal{G}|\geq|X_i|(\delta(G)+1)-|X_i|=|X_i|\delta(G)$ for any
$i\in\{2,\ldots,k\}$. Analogously, $|[X_i,\bar{X_2}]_\mathcal{G}|\geq|X_i|\delta(G)$, for any $i\in\{1,\ldots,k\}\setminus\{2\}$. Denote $M=\sum_{i=3}^k|X_i|$. Then $M\geq k-2$. By lemma 2.2, we have
$$
\begin{aligned}
\sum_{i=1}^{k}\sum_{j\in\{1,\ldots,k\}\setminus\{i\}}|[X_i,\bar{X_{j}}]_\mathcal{G}|
&= \sum_{i=2}^k|[X_i,\bar{X_1}]_\mathcal{G}|+\sum_{i\in\{1,\ldots,k\}\setminus\{2\}}|[X_i,\bar{X_2}]_\mathcal{G}|
+\sum_{i=3}^{k}\sum_{j\in\{3,\ldots,k\}\setminus\{i\}}|[X_i,\bar{X_{j}}]_\mathcal{G}|\\
&\geq(|X_2|+\sum_{i=3}^k|X_i|)\delta(G)+(|X_1|+\sum_{i=3}^k|X_i|)\delta(G)+\frac{(k-2)(k-3)}{2}2\lambda(G)\\
&\geq2(M+1)\delta(G)+(k-2)(k-3).
\end{aligned}
$$

Recall that $2\leq k\leq n-1$ and  $\delta(G)\geq 1$, then $Mn\geq Mk+M$ and $M\delta(G)\geq M$. If $M=0$, then $k=2$ and $V(D_1)=X_1\cup X_2$. Therefore, $|S|\geq\xi(\mathcal{G})=2n\delta(G)+2n-4$. Otherwise, assume $M\geq1$. Then by $(7)$, we have
$$
\begin{aligned}
|S|&\geq k(\delta(G)+1)+(M+2)[(\delta(G)+1)(n-k)-1]+2(M+1)\delta(G)+(k-2)(k-3)\\
&=(Mn-Mk+2M+2n-k+2)\delta(G)+k^2-6k-Mk+Mn-M+2n+4\\
&\geq(2n+2M)\delta(G)+2n+(k-5)(k-1)-1.\\
\end{aligned}
$$

Since $(k-5)(k-1)\geq-4$ for $k\geq 2$, we obtain $|S|\geq2n\delta(G)+2n+(k-5)(k-1)+2M\delta(G)-1>2n\delta(G)+2n-4$.

\noindent{\bf Subcase 2.2.3.} There is only one integer in $\{1,\ldots,k\}$, say $1$, such that 
$|X_1|=1$.

If there is a $|X_i|$ such that $|X_i|\geq3$, then $\sum_{i=1}^k|X_i|\geq2k$. By a similar argument as Subcase 2.2.1, we can also obtain $|S|>2n\delta(G)+2n-4$. Thus we assume $|X_i|=2$ for $2\leq i\leq k$.

If $\delta(G)=1$, then $\lambda(G)=\delta(G)=1$. By lemma 2.2, we have
$$
\begin{aligned}
\sum_{i=1}^{k}\sum_{j\in\{1,\ldots,k\}\setminus\{i\}}[X_i,\bar{X_{j}}]_\mathcal{G}|
&\geq\sum_{i=2}^{k}\sum_{j\in\{2,\ldots,k\}\setminus\{i\}}|[X_i,\bar{X_{j}}]_\mathcal{G}|
+\sum_{i=2}^k|[X_i,\bar{X_1}]_\mathcal{G}|\\
&\geq\frac{(k-1)(k-2)}{2}2\lambda(G)+(k-1)\cdot2\delta(G)\\
&=k(k-1).
\end{aligned}
$$

Since  $\sum_{i=1}^k|X_i|=2k-1$, by $(6)$, we have
$$
\begin{aligned}
|S|&\geq k\lambda(G)+(2k-1)(\delta(G)+1)(n-k)+k(k-1)\\
&=-3k^2+2k+4kn-2n.\\
\end{aligned}
$$

Set $f(k)=-3k^2+2k+4kn-2n$.
Recall that $2\leq k\leq n-1$ and $n\geq4$. Then we obtain $|S|\geq min\{f(2),f(n-1)\}=min\{6n-8,(n+1)(n-3)+4n-2\}>4n-4=2n\delta(G)+2n-4$.

Now we consider $\delta(G)\geq2$. Since $|X_i|\leq2$ for any $i\in\{1,\ldots,k\}$, we obtain $|[X_i,\bar{X_j}]|+|[X_j,\bar{X_i}]|\geq2\delta(G)$ for $1\leq i\neq j\leq k$. Thus
$$
\begin{aligned}
\sum_{i=1}^{k}\sum_{j\in\{1,\ldots,k\}\setminus\{i\}}|[X_i,\bar{X_{j}}]_\mathcal{G}|
\geq\frac{k(k-1)}{2}2\delta(G)=k(k-1)\delta(G).
\end{aligned}
$$

By using inequality $(7)$, we obtain
$$
\begin{aligned}
|S|&\geq k(\delta(G)+1)+(2k-1)[(\delta(G)+1)(n-k)-1]+k(k-1)\delta(G)\\
&> k(\delta(G)+1)+k[(\delta(G)+1)(n-k)-1]+k(k-1)\delta(G)\\
&=kn\delta(G)+kn-k^2\\
&=(k-2)(n\delta(G)+n-(k+2))+2n\delta(G)+2n-4\\
&\geq2n\delta(G)+2n-4.
\end{aligned}
$$
This proof is thus complete.
\end{proof}

Since $\xi(G\boxtimes K_n)=2n\delta(G)+2n-4$, Theorem 3.7 implies the following corollary.

\begin{cor}\label{Corollary3.8}
Let $G$ be a connected nontrivial graph with $m$ vertices. If min$\{n^2\lambda, (n-1)(m+2e(G))\}\geq 2n\delta(G)+2n-4$, then $G\boxtimes K_n$ is maximally restricted edge-connected, where $n\geq4$.
\end{cor}

By checking through the proof of Theorem 3.7, we find that $|S|\geq$min$\{n^2\lambda, (n-1)(m+2e(G))\}$ or 
$|S|>2n\delta(G)+2n-4$ when both $|V(D_1)|\geq3$ and $|V(D_2)|\geq3$. Thus we have the following corollary.

\begin{cor}\label{Corollary3.9}
Let $G$ be a connected nontrivial graph with $m$ vertices. If min$\{n^2\lambda, (n-1)(m+2e(G))\}> 2n\delta(G)+2n-4$, then $G\boxtimes K_n$ is super restricted edge-connected, where $n\geq4$.
\end{cor}

\end{document}